# IEEE Copyright Notice

# Stability Region Estimation Under Low Voltage Ride Through Constraints using Sum of Squares

Chetan Mishra, James S. Thorp, Virgilio A. Centeno
Bradley Dept. of Electrical and Computer Engineering,
Virginia Polytechnic Institute and State University,
Blacksburg-24061, Virginia, USA
Email: {chetan31, jsthorp, virgilio}@vt.edu

Anamitra Pal
School of Electrical, Computer, and Energy Engineering,
Arizona State University,
Tempe-85287, Arizona, USA
Email: Anamitra.Pal@asu.edu

*Abstract*— **The increasing penetration of inverter based renewable generation (RG) in the form of solar photo-voltaic (PV) or wind has introduced numerous operational challenges and uncertainties. According to the standards [1], [2], these generators are made to trip offline if their operating requirements are not met. In an RG-rich system, this might alter the system dynamics and/or cause shifting of the equilibrium points to the extent that a cascaded tripping scenario is manifested. The present work attempts at avoiding such scenarios by estimating the constrained stability region (CSR) inside which the system must operate using maximal level set of a Lyapunov function estimated through sum of squares (SOS) technique. A time-independent conservative approximation of the LVRT constraint is initially derived for a classical model of the power system. The proposed approach is eventually validated by evaluating the stability of a 3 machine test system with trip-able RG.**

*Keywords— Constrained stability region (CSR), direct methods, low voltage ride through (LVRT), sum of squares (SOS).*

## I. Introduction

Renewable generation (RG) poses numerous reliability challenges to successful power system operation. Besides the uncertainty in output, which can be dealt with up to a certain extent with forecasting, it is the response of inverter protection to system disturbances that has emerged as a major threat to system stability [3]. Traditionally inverter based generation was not expected to support the grid during disturbances and thus was tripped offline when such a scenario arose. However, with increase in RG penetration, this tripping could translate to losing large amounts of generation, more than what the system might be capable of sustaining, at the time of disturbance. A major development that occurred to safeguard against this problem was the development of inverter ride-through features which, as the name suggests, made the inverters stay online (ride through) during disturbances. The ride through standards would clearly define operating requirements (in terms of voltage and frequency) for staying online [1], [2].

Direct methods are one of the most popular methods used in industry and academia alike for analyzing transient stability of power systems [4]. However, direct methods are traditionally focused on capturing the mode that resulted in loss of synchronism of *synchronous generators*. Furthermore, their success heavily depends on the numerous system related assumptions that are made [5]. Loss of synchronism is characterized by a system trajectory leaving the stability region (SR) of the stable equilibrium point (SEP) of interest. Another way a system can leave the SR is when it transitions to another system having a different SEP. The latter can happen in a self-protecting system, such as a power system that is operating with considerable RG penetration as these generators are *constrained* to trip if the ride-through standards are not met during the disturbance. While at first it may not occur to be a threat, this constrained tripping in an RG-rich system could initiate a potentially cascaded tripping scenario that can have disastrous consequences. Making transient stability assessments for such constrained systems is also extremely difficult as along with loss of synchronism of synchronous generators, the direct methods must also estimate tripping of renewable generators.

The only previous research on this topic [6] converted such constrained systems into unconstrained systems, paving the way for energy function based methods like boundary controlling unstable (BCU) equilibrium point method [7] to become applicable. Also proposed was a time independent formulation of the low voltage ride through (LVRT) constraint. In this paper, we present an alternate formulation for solving this problem; one that is based on a Lyapunov function approach [8]. We will demonstrate the application of sum of squares (SOS) techniques in estimating the CSR for LVRT constrained systems.

The rest of the paper is structured as follows. The idea of CSR along with the modeling of LVRT constraint for a classical power system model is presented in Section II. The formulation of estimation of constrained region of attraction (ROA) in the form of SOS optimization is proposed in Section III. Lastly, effectiveness of the proposed methodology is demonstrated in Section IV and validated against time-domain simulations.

## II. Low Voltage Ride Through (LVRT) Constraint

### A. Constrained systems and stability region

Equality constraints when present, force the systems to evolve over a manifold. An example is the power system which is constrained to satisfy balance of nodal injections. Systems with equality constraints are usually modeled using differential algebraic equations. It is also common to find systems having inequality constraints. Such constraints can result from physical limitations of the system under study, requirement to meet some operating standard, or due to a preference of the system designer. An example is the condition that no RG should trip during a low voltage situation, which can be mathematically expressed as,

$$|v(t)| > LVRT(t) \qquad (1)$$

In (1), $v(t)$ is the voltage at the point of interconnection for the given inverter based generation and $LVRT(t)$ is the voltage reference value from the corresponding LVRT curve. When dealing with constrained systems, it is important to discuss the idea of feasible region of state space which is a collection of all points that do not violate any of the constraints. In order to do so, consider a generalized autonomous constrained system with state vector $x$ governed by the following equations:

$$\dot{x} = f(x) \qquad (2)$$
$$0 = g(x) \qquad (3)$$
$$0 \leq h(x) \qquad (4)$$

It should be noted that (2)-(4) are set-up such that $\forall x \text{ subject to } g(x) = 0$, $\frac{\partial g(x)}{\partial x}f(x) = 0$ i.e. if the system initializes on the manifold defined by $g(x) = 0$, then it can only evolve over that manifold. On the other hand, the inequality constraint, $h(x) \geq 0$, needs to be satisfied, but does not influence the system trajectories. This implies that the equilibrium points for the unconstrained system that lie inside the feasibility region make up the equilibrium points of the constrained system. Given a stable equilibrium point $x_s$ for this system, the CSR is defined as a collection of all those points in state space from which the emerging trajectories remain totally inside the feasibility region and eventually converge to $x_s$.

We estimate the CSR using the largest invariant set containing $x_s$ which does not intersect the infeasible part of state space. A Lyapunov function $V(x)$ [8], is estimated such that its corresponding level set $\{x|V(x) \leq 1\}$ serves as the invariant set to be found. We are only concerned with the portion of the state space $D = \{x \mid g(x) = 0\}$ where the trajectories are constrained to evolve. Assuming $x_s = 0$, the Lyapunov conditions to be fulfilled are [9]:

1. $V(x) > 0 \ \forall x \in D\setminus\{0\}$ and $V(0) = 0$
2. $\dot{V}(x) = \frac{\partial V}{\partial x}f(x) < 0 \ \forall \ x \in D\setminus\{0\} \cap \{x|V(x) \leq 1\}$
3. $h(x) \geq 0 \ \forall \ x \in D \cap \{x|V(x) \leq 1\}$

*B. LVRT constraint modeling for the classical power system*

LVRT curve is a time-dependent lower limit of the voltage at the renewable generator bus. Now, having time-independent constraints is much simpler to deal with. We will be taking a conservative approach in this study by approximating LVRT curve as a constant value equal to its maximum. While this constraint might appear to be very limiting, for systems with static load models the voltage recovery at the instance of fault clearing is almost instantaneous. Therefore, this will not be a conservative approximation for such systems. Mathematically, for the $j^{th}$ inverter based generation this can be written as,

$$|v_j(t)| > Max_t \left(LVRT_j(t)\right) \ \forall t \qquad (5)$$

We will first use the network reduced model of the classical power system. In accordance with the common modeling practice, inverter based generation will be modeled as a negative real load. Thus, inverter dynamics will not be modeled. Since we will be using a sum of squares (SOS) based approach which requires polynomial systems, the first task is to convert the system into a polynomial system through variable transformation. We will use the transformation described in Table I [10], where superscript $s$ indicates the SEP of interest and $n_g$ refers to the synchronous machine.

TABLE I. VARIABLE TRANSFORMATION [10]

| New Variable | Function of Original States |
|---|---|
| $z_i$ | $\omega_{in_g} - \omega_{in_g}^s(0) = \omega_{in_g}$ |
| $z_{n_g+2i-2}$ | $\sin(\delta_{in_g} - \delta_{in_g}^s)$ |
| $z_{n_g+2i-1}$ | $1 - \cos(\delta_{in_g} - \delta_{in_g}^s)$ |

The classical power system model for $n_g$ in a single machine reference frame for uniform damping is shown below.

$$\dot{\delta}_{in_g} = \omega_{in_g} \qquad (6)$$

$$M_i \dot{\omega}_{in_g} = P_{m_i} - \left(\sum_{j=1:n_g} E_i E_j \left(G_{ij} \cos\left(\delta_{in_g} - \delta_{jn_g}\right) + B_{ij} \sin\left(\delta_{in_g} - \delta_{jn_g}\right)\right)\right) + -\frac{M_i}{M_{n_g}}(P_{m_{n_g}} - \sum_{j=1:n_g} E_{n_g} E_j \left(G_{n_g j} \cos\left(-\delta_{jn_g}\right) + B_{n_g j} \sin\left(-\delta_{jn_g}\right)\right)) - D_i \omega_{in_g} \qquad (7)$$

$$\delta_{in_g} = \delta_i - \delta_{n_g}, \omega_{in_g} = \omega_i - \omega_{n_g} \qquad (8)$$
$$i = 1,2\ldots(n_g-1)$$

Using the variable transformation described in Table I, the system becomes a constrained system with the origin being the relevant SEP as shown in (9)-(12).

$$\dot{z}_i = \sum_j a_{ij} z_{n_g+2i-2} z_{n_g+2j-1} + \qquad (9)$$
$$\sum_j b_{ij} z_{n_g+2i-2} z_{n_g+2j-2} +$$
$$\sum_j c_{ij} z_{n_g+2i-1} z_{n_g+2j-2} +$$
$$\sum_j d_{ij} z_{n_g+2i-1} z_{n_g+2j-1} +$$
$$\sum_j e_{ij} z_{n_g+2j-1} + \sum_j f_{ij} z_{n_g+2j-2}$$

$$\dot{z}_{n_g+2i-2} = (1 - z_{n_g+2i-1}) \times z_i \qquad (10)$$
$$\dot{z}_{n_g+2i-1} = z_{n_g+2i-2} \times z_i \qquad (11)$$
$$0 = z_{n_g+2i-2} + \left(1 - z_{n_g+2i-2}\right)^2 - 1 \qquad (12)$$
$$i = 1,2\ldots n_g - 1$$

Since the network buses are reduced, the expression for voltages at the inverter buses needs to be derived as an explicit function of states $(\delta, \omega)$. The LVRT constraint also needs to be a polynomial function of the transformed states. Using the network reduction process described in [11], we get

$$\begin{bmatrix} 0 \\ I_g \end{bmatrix} = \begin{bmatrix} Y_{11} & Y_{12} \\ Y_{21} & Y_{22} \end{bmatrix} \begin{bmatrix} v \\ E \end{bmatrix} \qquad (13)$$

where $I_g$ is the injected current vector, $E$ is the internal EMF vector for generators (excluding inverter based generators which are modeled as negative loads), $v$ is the voltage vector for network buses (to be reduced), and $Y$'s are the components of the admittance matrix. Re-arranging (13), we get,

$$[v]_{n_{bus} \times 1} = \qquad (14)$$
$$-[Y_{11}]^{-1}[Y_{12}] \begin{bmatrix} |E_1|\left(\cos\left(\delta_1 - \delta_{n_g}\right) + j\sin\left(\delta_1 - \delta_{n_g}\right)\right) \\ \vdots \\ |E_{n_g}|\left(\cos\left(\delta_{n_g} - \delta_{n_g}\right) + j\sin\left(\delta_{n_g} - \delta_{n_g}\right)\right) \end{bmatrix}$$

Substituting the transformed states in (14), we get (15).
$$[v]_{n_{bus} \times 1} = [K]_{n_{bus} \times n_g} \times \begin{bmatrix} 1 - z_{n_g+1} + jz_{n_g+2} \\ 1 - z_{n_g+3} + jz_{n_g+4} \\ \vdots \end{bmatrix} \quad (15)$$

In (15), $K$ is a constant matrix consisting of terms from $E$, $Y_{11}^{-1}$ and $Y_{12}$. Now, as the bus voltage magnitudes involve a square root term they are not polynomial functions. Therefore, we use the LVRT constraint in terms of the square of bus voltage magnitudes which then become polynomial functions of the transformed variables. Let $n_{inv}$ inverter based generators with respective LVRT curves be present. Then, the LVRT constraint for each will be written in terms of lower limit on the voltage magnitude squared as shown below.

$$h_j(z_{ng}, z_{ng+1} \dots) = |v_j|^2 \quad (16)$$
$$- Max_t \left( LVRT_j(t) \right)^2$$
$$\geq 0 \; \forall j \in [1, n_{inv}]$$

III. CONSTRAINED ROA ESTIMATION USING SUM OF SQUARES

Sum of squares (SOS) based SR estimation has shown promising results in comparison to pre-existing Lyapunov based techniques for power systems [10]. Owing to the ease with which complex systems can be handled in the SOS approach, we will be using this technique for estimating our CSR. Before doing so, we provide a brief background on SOS.

A. *Relevant concepts in sum of squares (SOS)*

A function $F(x)$ is called a sum of squares if it can be written as $F(x) = \sum f_i(x)^2$. For a polynomial $F(x)$, the SOS check is equivalent to finding a symmetric positive semidefinite matrix $Q$(gram matrix) such that $F(x) = z^T(x)Qz(x)$ where $z(x)$ are all the monomials of degree $\leq$ half the degree of $F(x)$ [12]. It has been shown that $Q$ belongs to an affine subspace of symmetric matrices and is written as

$$Q = Q_0 + \sum_{i=1}^{n} \lambda_i Q_i \quad (17)$$

where $z(x)^T Q_0 z(x) = F(x)$ and $Q_i$ are basis matrices fulfilling the condition $z^T Q_i z = 0$. Thus, SOS problem can be converted to a linear matrix inequality (LMI) problem of finding $\lambda's$.

$$Q_0 + \sum_{i=1}^{n} \lambda_i Q_i \geq 0$$

Now, Positivestellensatz (P-Satz) Theorem [13] states that given polynomials $\{f_1(x), f_2(x) \dots f_s(x)\}$, $\{g_1(x), g_2(x) \dots g_t(x)\}$, $\{h_1(x), h_2(x) \dots h_u(x)\}$ the following are equivalent:

$$\left\{ x \in R^n \begin{vmatrix} f_1(x) \geq 0, f_2(x) \geq 0 \dots f_s(x) \geq 0 \\ g_1(x) \neq 0, g_2(x) \neq 0 \dots g_t(x) \neq 0 \\ h_1(x) = 0, h_2(x) = 0 \dots h_u(x) = 0 \end{vmatrix} \right\} \quad (18)$$
$$= \emptyset$$

Polynomials $F \in Cone(f_1(x), f_2(x) \dots f_s(x))$, $G \in Monoids(g_1(x), g_2(x) \dots g_t(x))$, $H \in Identity(h_1(x), h_2(x) \dots h_u(x))$ exist such that
$$F + G^2 + H = 0 \quad (19)$$

This theorem is a powerful tool for modeling the constraints as SOS constraints. For the rest of the paper, we will be using $z$ and not $x$ to represent the state vector.

B. *Local estimate of Lyapunov function*

To get an initial Lyapunov function and its corresponding ROA estimate, we use an iteration of expanding D algorithm [14]. We start with a known polynomial $p(z)$ and a semi-algebraic set $P_\beta = \{z|p(z) \leq \beta\}$ and expand it such that $\forall z \neq 0 \in P_\beta$, unknown function $V(z)$ is decreasing $\frac{\partial V}{\partial z}\dot{z} < 0, V > 0$, and $h(z) \geq 0$. Here it should be kept in mind that since the system dynamics are constrained to the manifold $A = \{z|g(z) = 0\}$, the constraints are only required to be fulfilled on it. In order to apply P-Satz theorem, we represent the constraints as set emptiness conditions as shown below.

$$\max \beta$$
$$st.$$
**Constr:** $V(z)$ cannot be $\leq 0 \; \forall z \neq 0 \in P_\beta \cap A$:
**P-Satz:** $F_1 \in Cone(-V, \beta - p(z)) + G_1^2 \in Monoid(l_1(z)) + H_1 \in Identity(g(z)) = 0$
**Constr:** $\dot{V}(z)$ cannot be $\geq 0 \; \forall z \neq 0 \in P_\beta \cap A$:
**P-Satz:** $F_2 \in Cone(\dot{V}, \beta - p(z)) + G_2^2 \in Monoid(l_2(z)) + H_2 \in Identity(g(z)) = 0$
$$-s_6(\beta - p) - \dot{V} - \lambda_2^T g - l_2 = \Sigma_n$$
**Constr:** $h_j(z)$ cannot be $< 0 \; \forall z \in P_\beta \cap A$:
**P-Satz:** $F_3 \in Cone(-h_j(z), \beta - p(z)) + G_3^2 \in Monoid(-h_j(z)) + H_2 \in Identity(g(z)) = 0$
$$\forall j \in [1, n_{inv}]$$

Due to the presence of $z \neq 0$ term in the first two constraints, those are not semi-algebraic. This is converted to $l_k(z) \neq 0$ where $l_k(z) = \epsilon z^T z$, $\epsilon$ is a small number [14]. The modified SOS constraint conditions are given by,

$$\max \beta$$
$$V, \dot{V}, s_2, s_6, s_{9j}, \lambda_1, \lambda_2, \lambda_3$$
$$-s_2(\beta - p) + V - \lambda_1^T g - l_1 \text{ is SOS}$$
$$-s_6(\beta - p) - \dot{V} - \lambda_2^T g - l_2 \text{ is SOS}$$
$$-s_{9j}(\beta - p) + h_j - \lambda_{3j}^T g \text{ is SOS}$$
$$s_i's \text{ are SOS}$$
$$\forall j \in [1, n_{inv}]$$

The functions $F, G$ & $H$ in formulation according to the P-Satz theorem are not unique. Also, there can be terms which are a multiplication of multiple unknown functions. For instance, in the above expression, besides an unknown $V$, there are other unknown multiplier functions ($\lambda's$ and $s's$). These cannot be solved through SOS programming as the problem cannot be converted to a corresponding LMI. Therefore, the degrees need to be chosen carefully to make it possible to find a feasible solution [14]. For example, for a known function $V$, if the constraints are $s_1 - V$ is SOS and $s_1$ is SOS, first constraint can only be an SOS if its degree is even. So, if $V$ has an odd degree, degree of $s_1 >$ degree of $V$. As such, some of the rules we follow when choosing the degrees for multiplier functions are:

a. To limit the size of the problem, each term in the constraint expression is chosen to have the same/similar degree with the overall degree of the known positive terms $\geq$ degrees of other terms.
b. Degrees of SOS multipliers are even by definition and so each constraint's degree is even as well.
c. Since $V$ is 0 at $z = 0$, $V$ does not have constant terms.

For further reading on simplifications, please see [12].

*C. Iterative ROA estimate using expanding interior algorithm*

From Section III.B, we get an initial Lyapunov function $V$. One way of estimating ROA using a given Lyapunov function is through its level sets. Ref. [14] states that if $D \subseteq \mathbb{R}^n$ be a domain containing equilibrium $x = 0$ of the system $\dot{x} = f(x)$ with $V(x)$ being the corresponding Lyapunov function defined in $D$, then any region $\Omega_\beta := \{x \in D | V(x) \leq \beta\}$ is a positive invariant region contained in the equilibrium's ROA. Thus, for a given Lyapunov function, its largest level set contained inside the region $D$ gives an estimate of its ROA. Since our system also has inequality constraints, another condition to be fulfilled by the level set is that it stays within the feasible region.

We use the expanding interior algorithm [14] to convert the problem of finding the biggest level set for an unknown Lyapunov function into an SOS optimization problem. The idea being that a semi-algebraic set $P_\beta$ contained inside a given local estimate for ROA can be expanded. During the expansion process, the Lyapunov function as well as its associated ROA estimate are allowed to change as long as they totally contain $P_\beta$. Now, this makes the problem heavily dependent on the shape of $P_\beta$ chosen. An alternative to this dependency was proposed in [10] where an outer iteration loop was added to change $P_\beta$ into the obtained Lyapunov level set at the end of the inner expansion process. The optimization problem with set emptiness ($\emptyset$) constraints after incorporating all the previously discussed conditions can be written as: Given $P_\beta = \{z | p(z) \leq \beta\}$ with radius $\beta$ and an unknown Lyapunov function $V$ with the associated SR estimate given by $\{z | V(z) \leq 1\}$,

$$\text{max. } \beta$$
$$\text{s.t.}$$

$V(z)$ cannot be $\leq 0$ for $z \neq 0$ : $\{V(z) \leq 0, g(z) = 0, z \neq 0\} = \emptyset$ (20)

$P_\beta$ Is contained inside the SR estimate $(V(z) \leq 1)$ : $\{p(z) \leq \beta, g(z) = 0, V(z) \geq 1, V(z) \neq 1\} = \emptyset$ (21)

Inside the SR estimate, $V(z)$ strictly decreases along all trajectories: (22)
$\{V(z) \leq 1, \dot{V}(z) \geq 0, g(z) = 0, z \neq 0\} = \emptyset$

Inside the SR estimate, $h_j(z)$ is never $< 0$: (23)
$\{V(z) \leq 1, -h_j(z) \geq 0, -h_j(z) \neq 0, g(z) = 0\} = \emptyset$
$\forall j \in [1, n_{inv}]$

Using P-Satz theorem as well as the definitions for Cones, Multiplicative Monoids and Identity generator, we can convert these constraints into the following equivalent conditions:

$$s_1 - s_2 V + \lambda_1^T g + l_1^{2k_1} = 0 \tag{24}$$
$$s_3 + s_4(\beta - p) + s_5(V - 1) \tag{25}$$
$$+ s_6(\beta - p)(V - 1) + \lambda_2^T g$$
$$+ (V - 1)^{2k_2} = 0$$
$$s_7 + s_8(1 - V) + s_9 \dot{V} + s_{10}(1 - V)\dot{V} + \lambda_3^T g \tag{26}$$
$$+ l_2^{2k_3} = 0$$
$$s_{11j} + s_{12j}(1 - V) - s_{hj}(1 - V)h_j - s_{14j} h_j \tag{27}$$
$$+ \lambda_{hj}^T g + (-h_j)^{2k} = 0$$

Here, we deliberately leave out the higher degree terms in all the constraints (for example $V^2$ in the first constraint) to make LMI formulation possible as discussed before. Since a product of SOS functions is also SOS, in the first and third constraints, assuming $k_1 = k_2 = 1$ and $s_{10} = 0, s_1 = s_1 \times l_1, s_2 = s_2 \times l_1, \lambda_1 = \lambda_1 \times l_1, s_7 = s_7 \times l_2, s_8 = s_8 \times l_2, s_9 = s_9 \times l_2, \lambda_3 = \lambda_3 \times l_2$, we will eliminate $l_1$ and $l_2$. Similarly in the second constraint, assuming $k_2 = 1, \lambda_2 = \lambda_2 \times (V - 1), s_3 = s_4 = 0$, we can take $(V - 1)$ common from each term. Finally, for the fourth constraint we assume $s_{11j} = s_{12j} = 0, \lambda_{hj} = \lambda_{hj} \times h_j$ and take $h_j$ common. Taken together, we get,

$$s_2 V - \lambda_1^T g - l_1 = s_1 = SOS \tag{28}$$
$$-s_6(\beta - p) - \lambda_2^T g - (V - 1) = s_5 = SOS \tag{29}$$
$$-s_8(1 - V) - s_9 \dot{V} - \lambda_3^T g = s_7 = SOS \tag{30}$$
$$-s_{hj}(1 - V) - \lambda_{hj}^T g + h_j = s_{14j} = SOS \tag{31}$$

Here, it can be seen that some of the terms have 2 unknown functions multiplied (for example $s_9 \times \dot{V}$ in third constraint) which makes the problem unsolvable by semi-definite programming. To fix this problem, $\beta$ is maximized for a given fixed $p$ using the algorithm summarized below.

---

**Expanding Interior (fixed p, maximize $\beta$)**

i. Input $p, \beta^{(0)}$. Initialize $i = 0, V^{(0)}, \dot{V}^{(0)}$

ii. Substitute $V = V^{(i)}, \dot{V} = \dot{V}^{(i)}$ in constraints in Eqn. 28-31, line search over $\beta$ starting from $\beta^{(i)}$ till the SOS problem is infeasible. Let the final obtained values be $\beta', s_8', s_9', s_{hj}'$. Set $\beta^{(i+1)} = \beta'$

iii. Substitute $s_8 = s_8', s_9 = s_9', s_{hj} = s_{hj}'$ in constraints in Eqn 28-31, line search over $\beta$ starting from $\beta^{(i+1)}$ till the SOS problem is infeasible. Let the final obtained values be $\beta', V'$. Update $\beta^{(i+1)} = \beta', V^{(i+1)} = V', \dot{V}^{(i+1)} = \dot{V}'$.

iv. If $|\beta^{(i)} - \beta^{(i+1)}| > tolerance, i = i + 1$, Goto ii.

v. Stop.

---

As mentioned before, since the technique is dependent on the choice of the function $p$, a workaround is to update $p = V$ at the end of each expanding interior loop [10]. This gives the following overall algorithm.

Overall Algorithm

i. Initialize $j = 0, p^{(0)} = z^T z, \beta^{(0)} = 0$

ii. Run **Expanding Interior** with $p^{(j)}, \beta^{(j)}$ as inputs to get $V$.

iii. If $V \neq p^{(j)}$ then $p^{(j+1)} = V, \beta^{(j+1)} = 1, j = j + 1$, Goto ii.

iv. Stop

IV. RESULTS AND DISCUSSION

We will demonstrate the implementation as well as the effectiveness of the proposed technique on a standard 3 machine system [15] with RG connected as shown in Fig. 1. A small PV modeled as negative real load $(-0.4 + 0j)$ p.u. is added to bus 1 with the generator's output at that bus reduced accordingly. IEEE 1547 standard LVRT curve is used for the PV which has a maximum value of 0.85 p.u. (Fig. 1) Thus, the feasibility region is given by $v_1 \geq 0.85$ where $v_1$ is the voltage at bus 1. The post-fault system with the PV connected has the SEP at (0.3165,0.3451,0,0) resulting in the transformed state variables given in Table II. The chosen degrees for the multiplier functions are given in Table III.

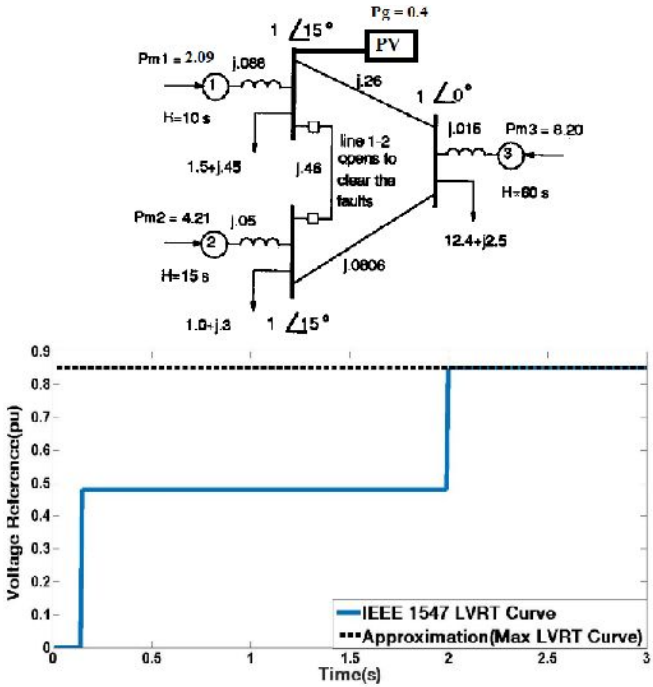

**Figure 1a. Three machine system with PV connected at bus 1; Figure 1b. IEEE 1547 standard LVRT curve**

The obtained estimated Lyapunov function for the constrained system is given by,

$$V(z) = 0.019 \times z_1^2 - 3.7 \times 10^{-3} \times z_1 \times z_2 + 0.025 \times z_1 \times z_3 - 0.019 \times z_1 \times z_4 - 0.018 \times z_1 \times z_5 + 0.012 \times z_1 \times z_6 + 3.3 \times 10^{-3} \times z_2^2 + 2.4 \times 10^{-3} \times z_2 \times z_3 - 4.4 \times 10^{-3} \times z_2 \times z_4 + 0.021 \times z_2 \times z_5 + 4.8 \times 10^{-3} \times z_2 \times z_6 + 1.0 \times z_3^2 + 0.21 \times z_3 \times z_4 - 0.14 \times z_3 \times z_5 - 4.0 \times 10^{-3} \times z_3 \times z_6 + 0.65 \times z_4^2 - 0.018 \times z_4 \times z_5 - 0.19 \times z_4 \times z_6 + 0.21 \times z_4 + 0.74 \times z_5^2 - 0.024 \times z_5 \times z_6 + 0.7 \times z_6^2 - 0.73 \times z_6$$

TABLE II. PROPOSED VARIABLE TRANSFORMATION

| New Variable | In terms of Original System State |
|---|---|
| $z_1$ | $\omega_{1n_g}$ |
| $z_2$ | $\omega_{2n_g}$ |
| $z_3$ | $\sin(\delta_{1ng} - 0.3165)$ |
| $z_4$ | $1 - \cos(\delta_{1ng} - 0.3165)$ |
| $z_5$ | $\sin(\delta_{2ng} - 0.3451)$ |
| $z_6$ | $1 - \cos(\delta_{2ng} - 0.3451)$ |

TABLE III. MULTIPLIER FUNCTION DEGREES

| Function | Degree | Function | Degree |
|---|---|---|---|
| $s_2$ | 0 | $\lambda_{21}, \lambda_{22}$ | 0 |
| $s_6$ | 0 | $\lambda_{31}, \lambda_{32}$ | 2 |
| $s_8$ | 2 | $s_{h1}$ | 0 |
| $s_9$ | 0 | $\lambda_{h_{11}}, \lambda_{h_{12}}$ | 0 |
| $\lambda_{11}, \lambda_{12}$ | 0 | | |

The estimated SRs $\{z|V(z) \leq 1\}$ in terms of original states for the given constrained system vs. the same system with no LVRT constraint are shown for $\omega = 2$ in Fig. 2. It can be noticed that since the feasibility region is more binding in the horizontal direction, the estimate hits that feasibility boundary first and stops expanding in that direction.

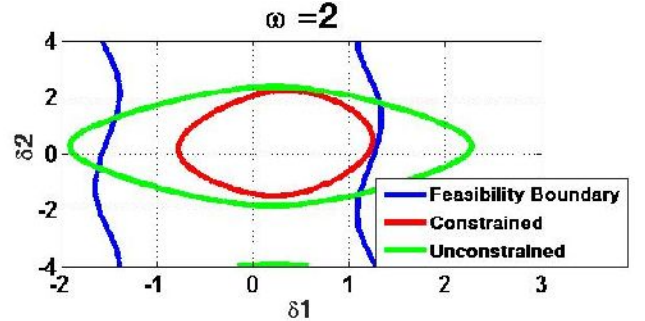

**Figure 2. Constrained vs Unconstrained ROA**

In order to validate the estimated constrained ROA, voltage at bus 1 was plotted for post fault trajectories starting from random points within it. As seen in Fig. 3 below, the voltages for all the trajectories never fall below 0.85 p.u. Also, the projection of trajectories emerging from the SR boundary on the angle plane shows that they are stable as well (Fig. 4). Now, in order to analyze the level of conservativeness of the proposed approach, we obtain the true CSR through time domain simulation using the maximum of LVRT curve. As discussed before, if the trajectory emerging from a point in state space intersects the infeasible region or does not return to the SEP, it cannot belong to the CSR. It can be observed from Fig. 5 that while the estimate covers a large portion of the one obtained through time domain simulation, since the SEP was closer to the right side of the feasibility boundary, when expanding the ROA estimate with SEP as the center, it stops as soon as it intersects the right side boundary. The proposed technique could yield conservative results for systems having slower voltage recovery, as it could be difficult to maintain the immediate post fault system voltage above the maximum value of LVRT curve. However, systems with static load models can be dealt with effectively using the proposed approach.

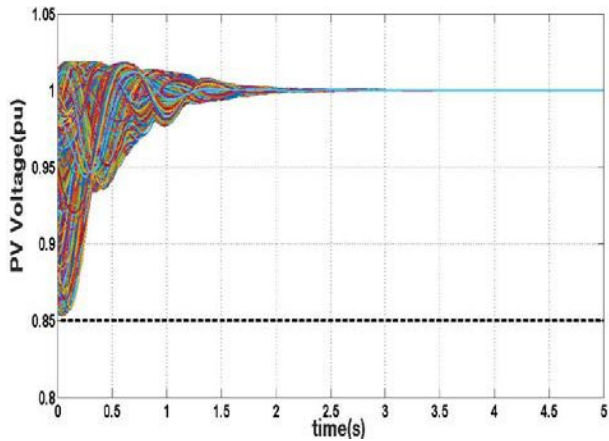

**Figure 3. Voltage along trajectories emerging from the estimated stability boundary**

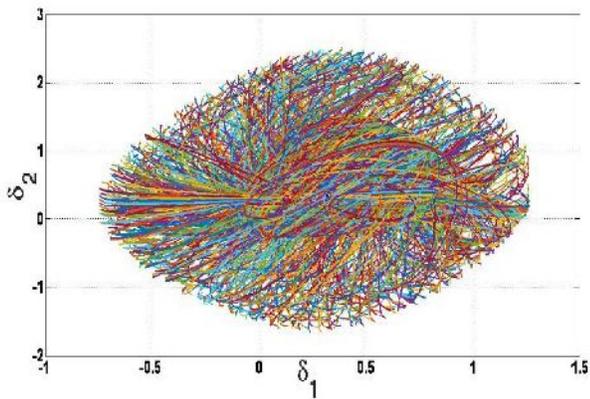

**Figure 4. Projection of trajectories from the estimated stability boundary on the angle plane**

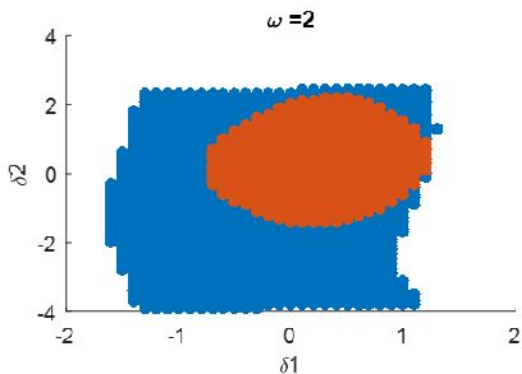

**Figure 5. Estimated constrained ROA (brown) vs actual (blue)**

## V. CONCLUSIONS AND FUTURE WORK

In this paper, estimation of stability region (SR) under the low voltage ride through (LVRT) constraint is tackled using a Lyapunov-based approach through SOS programming. It was seen that SOS serves as a powerful tool to deal with the analysis of complex dynamical systems. A conservative formulation of LVRT constraint followed by the derivation of SOS constraints was proposed for a classical power system model. It was seen from the results obtained for a 3 machine test system that the proposed approach gave satisfactory estimates of the CSR.

However, there is still scope for improvement. One major challenge is that due to computational limitations, it is extremely difficult to use the proposed approach in large scale systems. In this regard, the idea of decomposition using a vector Lyapunov function (found to be applicable to large interconnected systems in [16]) can be explored. With regards to the conservativeness in our approach, modeling of LVRT constraints similar to [6] could be tried. Lastly, since the final estimate depended heavily on the relative location of the SEP and the feasibility boundary, in the expanding interior algorithm, shifting of the center of the expanding region $p_\beta$ could be attempted.